\documentclass[11pt]{article}

\usepackage{amsmath,amssymb,amsthm}
\usepackage{mathrsfs}
\usepackage{booktabs}
\usepackage{xcolor}
\usepackage{graphicx}
\usepackage[margin=1in]{geometry}
\usepackage{enumitem}
\usepackage{hyperref}
\hypersetup{colorlinks=true,linkcolor=blue!60!black,citecolor=blue!60!black,urlcolor=blue!60!black}

\newcommand{\cA}{\mathcal{A}}
\newcommand{\cQ}{\mathcal{Q}}
\newcommand{\cH}{\mathcal{H}}
\newcommand{\cT}{\mathcal{T}}
\newcommand{\cP}{\mathscr{P}}
\newcommand{\cN}{\mathscr{N}}
\newcommand{\om}{\,o^-}

\newcommand{\N}{\mathbb{N}}

\newtheorem{theorem}{Theorem}[section]
\newtheorem{lemma}[theorem]{Lemma}
\newtheorem{proposition}[theorem]{Proposition}
\newtheorem{corollary}[theorem]{Corollary}
\newtheorem{definition}[theorem]{Definition}
\newtheorem{remark}[theorem]{Remark}
\newtheorem{conjecture}[theorem]{Conjecture}
\newtheorem{problem}[theorem]{Open Problem}

\title{Quaternary games with infinite mis\`ere quotients%
\thanks{The programs and results of this paper are available at
\url{https://github.com/liang-kai-ruqi/automata_and_games}.}}
\author{Kai Liang}
\date{\today}

\begin{document}
\maketitle

\begin{abstract}
We study the four wild quaternary games left open by the earlier work of
Plambeck and Siegel. For \texttt{0.3102} we discover a special regularity of
its two--heap case, give a complete solution of the two--heap case, and deduce
from it that its mis\`ere quotient is infinite. For \texttt{0.3122} and
\texttt{0.3123} we rigorously prove, by slice--wise algebraic periodicity and
hierarchical induction (computer--assisted), that their mis\`ere quotients are
infinite. For \texttt{0.3312}, as well as many other wild quaternary games, we
carry out computations and observations.
\end{abstract}



\section{Introduction}\label{sec:intro}

In normal play, the Sprague--Grundy theory assigns to every impartial game an
SG value (nim value); the SG value of a disjunctive sum is the nim--sum (bitwise
xor) of the SG values of the summands, and a position is a $\cP$--position if and
only if its SG value is $0$. For an octal game this means that the outcomes of
all positions are determined by the SG values of the single--heap positions.
Although many cases of octal games remain open, their special case without
splitting, namely the quaternary games, always have eventually periodic
single--heap positions, so that Sprague--Grundy theory yields a complete
solution.

In mis\`ere play, however, Sprague--Grundy theory fails: in general no single
``value'' with the above property exists, and even quaternary games become
intricate. For this reason Plambeck and Siegel introduced the theory of the
mis\`ere quotient \cite{PS,MisQuot2006}: two positions are equivalent if they
have the same outcome in every context, and the resulting quotient is a
commutative monoid $\cQ$. The mis\`ere quotient plays the role that the nim
group plays in normal play: the outcome of a disjunctive sum is a function of
the product (in $\cQ$) of the images of the summands.

Plambeck and Siegel's results \cite{PS,miseregames} show that almost all
quaternary games of the form
$0.\mathrm{xxxx}$ have a finite mis\`ere quotient and can therefore be solved
completely by computing the quotient; only four games remained open:
\[
\texttt{0.3102},\quad \texttt{0.3122},\quad \texttt{0.3123},\quad
\texttt{0.3312},
\]
and it was conjectured that their mis\`ere quotients are infinite.

This paper proves the infinitude of the mis\`ere quotient for three of them. For
\texttt{0.3102}, the simplest of the four,
we give a complete proof of the two--heap outcome: the outcome of $(x,y)$ is
$\cP$ if and only if $|y-x|\equiv 0,3\pmod 5$, and we deduce that its mis\`ere
quotient is infinite. For \texttt{0.3122} and \texttt{0.3123}, the two--heap
outcome is simply periodic, and the non--periodicity appears only on four--heap
slices; we rigorously prove that their quotients are infinite using slice--wise
algebraic periodicity and hierarchical induction (computer--assisted, see
Section~\ref{sec:proven}). The remaining game \texttt{0.3312} is still open: its
non--periodicity lies in the three--heap diagonal neighbourhood, where even the
diagonal reduction fails, so the above method does not apply. In addition, we
observe and summarize all quaternary games of code length at most $8$, and find
that as the code length grows, the number of games with a similarly intricate
structure grows rapidly.

\section{Preliminaries}\label{sec:prelim}

\subsection{Octal games and the mis\`ere outcome}\label{sec:outcome}

An \emph{octal game} $\Gamma$ is a take--and--break game on heaps of tokens: a
move on a single heap of size $n$ removes $k$ tokens and, according to a fixed
code, possibly splits the remaining $n-k$ tokens into several heaps. A
\emph{position} is a finite multiset of heaps; we represent it by a
non--increasing tuple of positive integers $(a_1,\dots,a_n)$ with
$a_1\ge a_2\ge\cdots\ge a_n$, and the empty tuple $()$ is the empty game. The
set of all positions is the free commutative monoid
\[
\cA=\bigl\{\,()\,\bigr\}\cup
\bigl\{\,(a_1,\dots,a_n): n\ge1,\ a_1\ge a_2\ge\cdots\ge a_n\ge1\,\bigr\},
\]
the \emph{heap algebra}, with multiplication the multiset union (written by
juxtaposition or $+$). We write $\cH=\{(n):n\ge 1\}$ for the set of single--heap
positions.

In mis\`ere play the last player to move loses. The \emph{mis\`ere outcome}
$\om(G)\in\{\cP,\cN\}$ of a position $G$ is $\cP$ (resp.\ $\cN$) if the
\emph{previous} (resp.\ \emph{next}) player can force a win. It satisfies the
recursion
\[
\om(G)=\cN\iff \text{some option } G' \text{ of } G \text{ has } \om(G')=\cP,
\]
with the convention $\om(())=\cN$ (the last move loses).

For a \emph{slice}---a four--heap position $(a,b,x,y)$ with two fixed heaps
$a,b$ and two variable heaps $x,y$---we write
\[
o(a,b,x,y):=\om(a,b,x,y)
\]
as a shorthand for $\om$ on slices; below $o$ and $\om$ are interchangeable on
slices (so $\om(2,2,h,h)$ and $o(2,2,h,h)$ denote the same object).

\subsection{The mis\`ere quotient}\label{sec:quotient}

Two positions $G,H\in\cA$ are \emph{equivalent}, written $G\equiv_{\cA}H$, if
\[
\om(G+X)=\om(H+X)\quad\text{for every } X\in\cA .
\]
If $G\not\equiv_{\cA}H$, a position $X$ with $\om(G+X)\neq\om(H+X)$ is called a
\emph{witness} of the non--equivalence of $G$ and $H$ (it ``witnesses'' that
$G$ and $H$ are not equivalent). This is a congruence on $\cA$. The
\emph{mis\`ere quotient} is
\[
\cQ=\cQ(\cA)=\cA/{\equiv_{\cA}},
\]
with $\cP$--portion
\[
\cP=\bigl\{[G]_{\equiv_{\cA}}:\om(G)=\cP\bigr\},
\]
and quotient map $\Phi:\cA\to\cQ$. Following \cite{PS} we define the
\emph{tame quotients}
\[
\cT_0=\cQ(0),\qquad \cT_n=\cQ({\ast2^{n-1}})\ \ (n\ge1),\qquad
\cT_\infty=\cQ(0,{\ast},{\ast2},{\ast3},\dots),
\]
where ${\ast k}$ denotes a Nim heap of size $k$; thus $\cT_n$ is the quotient of
the single--heap Nim position ${\ast2^{n-1}}$, and $\cT_\infty$ is the quotient
of mis\`ere Nim itself. A game is \emph{tame} if $\cQ$ is isomorphic to some
$\cT_n$ ($n\ge0$) or to $\cT_\infty$, and \emph{wild} otherwise.

\begin{theorem}[Periodicity Theorem, Fact 2.2 of \cite{PS}]\label{thm:periodicity}
If $\Phi((n+p))=\Phi((n))$ holds on an interval of length $p$, then it holds for
all sufficiently large $n$, and $\cQ$ is finite. Consequently $\cQ$ is infinite
if and only if $\Phi((n))$ is \emph{not} eventually periodic (since wildness
entails an infinite quotient, the single--heap images of a wild game are never
eventually periodic).
\end{theorem}

\subsection{Algebraic periodicity}\label{sec:algperdef}

The following definition is due to Liang \cite{gamedfa} (see also
\cite{MisQuot2006}).

\begin{definition}[Algebraic periodicity]\label{def:algper}
A game $\Gamma$ is \emph{algebraically periodic} if there exist natural numbers
$p_1,p_2\in\N$ and positive integers $t_1,t_2\in\N^+$ such that, for every
position $(a_1,a_2,\dots,a_n)\in\cA$ written with heaps in non--increasing order
($a_1\ge a_2\ge\cdots\ge a_n$) and with $a_1\ge p_1+t_1$, one of the following
holds:
\begin{enumerate}[label=(\roman*)]
\item \emph{single--heap reduction}: if $a_2<p_2+t_2$ (where $a_2=0$ when
$n=1$), then $\om(a_1,a_2,\dots,a_n)=\om(a_1-t_1,a_2,\dots,a_n)$;
\item \emph{two--heap reduction}: if $a_2\ge p_2+t_2$, then
$\om(a_1,a_2,\dots,a_n)=\om(a_1-t_2,a_2-t_2,a_3,\dots,a_n)$.
\end{enumerate}
The lexicographically smallest tuple $(p_1,t_1,p_2,t_2)$ is the
\emph{parameters}; $(p_1,t_1)$ and $(p_2,t_2)$ are the (pre--period, period) of
the two reductions respectively.
\end{definition}

Repeatedly applying the reductions reduces any position to one all of whose
heaps are smaller than $p_1+t_1$; the outcome of such ``corner'' positions is
given by a finite table. Thus an algebraically periodic game is ``solvable'' by
a finite amount of data.

\subsection{Quaternary games}\label{sec:games}

A \emph{quaternary game} has code digits in $\{0,1,2,3\}$: a move on a heap
either removes it (bit $0$) or replaces it by a single smaller heap (bit $1$);
there is no splitting (bit $2$). The four games we study are

\begin{center}
\begin{tabular}{cccc}
\toprule
Game & remove & subtract & single--heap period\\
\midrule
\texttt{0.3102} & $\{1,2\}$   & $\{1,4\}$     & $5$\\
\texttt{0.3122} & $\{1,2\}$   & $\{1,3,4\}$   & $7$\\
\texttt{0.3123} & $\{1,2,4\}$ & $\{1,3,4\}$   & $7$\\
\texttt{0.3312} & $\{1,2,3\}$ & $\{1,2,4\}$   & $3$\\
\bottomrule
\end{tabular}
\end{center}

Here ``remove $S$'' means a heap of size $s\in S$ may be removed (moved to
$0$), and ``subtract $T$'' means a heap of size $n>t$, $t\in T$, may be replaced
by a heap of size $n-t$. The fourth column is the period of the single--heap
outcome (computed directly).

\section{Rigorously proved results}\label{sec:proven}

\subsection{The two--heap outcome of \texttt{0.3102}}\label{sec:03102}

We give a complete proof of the two--heap outcome of \texttt{0.3102}.

\begin{theorem}\label{thm:03102}
In \texttt{0.3102}, for all $x,y\ge 3$,
\[
\om(x,y)=\cP \iff |y-x|\equiv 0,3 \pmod 5 .
\]
\end{theorem}

\begin{proof}
Write $f(x,y)=\om(x,y)$, identifying the empty game $0$ with a heap of size $0$
(so $f(0,y)=f(y,0)=\om(y)$). We first record two finite facts.

\textbf{(a) Single heaps.} A direct computation from the recursion gives
\[
\om(1)=\cP,\ \om(2)=\cN,\ \om(3)=\cP,\ \om(4)=\cN,\ \om(5)=\cN,
\]
\[
\om(6)=\cP,\ \om(7)=\cN,\ \om(8)=\cP,\ \om(9)=\cN,\ \om(10)=\cN.
\]
Indeed, each $\om(h)$ depends only on $\om(h-1)$, $\om(h-4)$, and the removal
condition $h\in\{1,2\}$; since the block $h=0,1,2,3,4$ coincides with
$h=5,6,7,8,9$, we obtain
\[
\om(h)=\cP \iff h\equiv 1,3 \pmod 5 \qquad(h\ge 1),
\]
with $\om(0)=\cN$.

\textbf{(b) Two heaps, one small.}
From the recursion for $f(1,y)$ and $f(2,y)$ one computes
\[
f(1,y)=\cP \iff y\equiv 2,4 \pmod 5 \qquad(y\ge 4),
\]
\[
f(2,y)=\cP \iff y\equiv 0 \pmod 5 \qquad(y\ge 5),
\]
with the small--$y$ exceptions
\[
f(1,0)=\cP,\quad f(1,1)=f(1,2)=f(1,3)=\cN,
\]
\[
f(2,0)=f(2,1)=\cN,\quad f(2,2)=\cP,\quad f(2,3)=f(2,4)=\cN .
\]

\textbf{(c) Induction.} We prove the theorem by strong induction on $s=x+y$.
The base case $s=6$ is the single position $(3,3)$: its only move (subtracting
$1$ from a heap) leads to $(2,3)$, and $f(2,3)=\cN$ by (b), so $f(3,3)=\cP$; and
$|3-3|=0\equiv 0\pmod 5$. Now assume the theorem for all $(x',y')$ with
$x',y'\ge 3$ and $x'+y'<s$, and consider $(x,y)$ with $x,y\ge 3$, $x+y=s$. By
symmetry we may assume $x\le y$, and put $d=y-x$. Since $x,y\ge 3$ there is no
removal; each move subtracts $1$ or $4$ from one heap, which changes $d$ by
$\pm 1$ modulo $5$, hence changes $d$ to $|d\pm 1|\pmod 5$.

\emph{Case 1: $d\equiv 0$ or $3\pmod 5$} (target $\cP$). Every successor has
$|d\pm 1|\equiv 1,2$, or $4\pmod 5$, hence is a $\cN$--position, provided the
successor has both coordinates $\ge 3$ (by induction) or, when a coordinate
drops to $2$ or $1$, by the table (b) (dropping to $2$ or $1$ can only come
from subtracting $1$ from $3$, or $4$ from $5$ or $6$):
\begin{itemize}
\item $x=3$ (subtracting $1$ from $x$ gives $(2,y)$): $d=y-3\equiv 0,3$ means
$y\equiv 3$ or $y\equiv 1\pmod 5$, so $y\not\equiv 0\pmod 5$ and $f(2,y)=\cN$;
\item $x=5$ (subtracting $4$ from $x$ gives $(1,y)$): $d=y-5\equiv 0,3$ means
$y\equiv 0$ or $y\equiv 3\pmod 5$, so $y\not\equiv 2,4\pmod 5$ and
$f(1,y)=\cN$;
\item $x=6$ (subtracting $4$ from $x$ gives $(2,y)$): $d=y-6\equiv 0,3$ means
$y\equiv 1$ or $y\equiv 4\pmod 5$, so $y\not\equiv 0$ and $f(2,y)=\cN$;
\item subtracting $4$ from $y$ so that $y=6$ drops to $2$ (here $x=3$, $d=3$):
this gives $(3,2)$, and $f(2,3)=\cN$ by (b).
\end{itemize}
(The remaining moves subtracting $1$ or $4$ from $y$---namely $y=3$, which is
the base case, and $y=5$ or $y=6$ with $x=6$, whose targets coincide with the
$x=5$ and $x=6$ cases above---produce no new target.)

\emph{Case 2: $d\equiv 1,2$, or $4\pmod 5$} (target $\cN$). There is a move to a
$\cP$--position:
\begin{itemize}
\item $d\equiv 1$ or $4$: subtract $1$ from the larger heap. The larger heap
is $\ge 4$ (since the smaller is $\ge 3$ and $d\ge 1$), so this is a legal move,
and it changes $d$ to $d-1\equiv 0$ or $3$, a $\cP$--position.
\item $d\equiv 2$: subtract $1$ from the smaller heap. This changes $d$ to
$d+1\equiv 3$. If the smaller heap is $>3$, the successor has both coordinates
$\ge 3$ and is $\cP$ by induction. If the smaller heap is $3$, the successor is
$(2,y)$ with $y=3+d\equiv 5\equiv 0\pmod 5$, so $f(2,y)=\cP$ by (b).
\end{itemize}

The induction is complete.
\end{proof}

\begin{corollary}\label{cor:03102}
In \texttt{0.3102}, for every $h\ge 0$,
\[
\om(h,h+2)=\cN,\qquad \om(h+5,h+2)=\cP .
\]
\end{corollary}

\subsection{The infinitude of the mis\`ere quotient of \texttt{0.3102}}\label{sec:slicecriterion}

\begin{lemma}\label{lem:bidirectional}
Let $\Gamma$ be an octal game whose mis\`ere quotient $\cQ$ is finite. Then there
exist $p,N$ such that for every fixed $a,b$ and every $x,y\ge N$,
\[
\om(a,b,x,y)=\om(a,b,x+p,y)=\om(a,b,x,y+p).
\]
\end{lemma}

\begin{proof}
By Theorem~\ref{thm:periodicity}, $\Phi((n+p))=\Phi((n))$ for all $n\ge N$ for
some $p,N$. Then $\om(a,b,x,y)=\cP$ iff
$\Phi(a)\Phi(b)\Phi(x)\Phi(y)\in\cP$, and $\Phi((x+p))=\Phi(x)$, so the product
is unchanged upon replacing $x$ by $x+p$; likewise for $y$.
\end{proof}

\begin{definition}\label{def:bidir}
A slice $o(a,b,x,y)$ is \emph{bi--directionally eventually periodic} if there
exist $p,N$ such that $o(a,b,x+p,y)=o(a,b,x,y+p)=o(a,b,x,y)$ for all $x,y\ge N$.
\end{definition}

By Lemma~\ref{lem:bidirectional}, finiteness of $\cQ$ forces every slice to be
bi--directionally eventually periodic. We show that the slice $(0,0)$ of
\texttt{0.3102} is not.

\begin{proposition}\label{prop:03102slice}
In \texttt{0.3102}, the slice $o(0,0,x,y)=\om(x,y)$ satisfies
(i) \emph{diagonal periodicity}: $o(0,0,x+t,y+t)=o(0,0,x,y)$ for every $t\ge 0$
and all $x,y\ge 3$; but
(ii) \emph{axis aperiodicity}: for every $t\ge 1$ there exist arbitrarily large
$x,y\ge 3$ with $o(0,0,x+t,y)\neq o(0,0,x,y)$.
\end{proposition}

\begin{proof}
By Theorem~\ref{thm:03102}, $o(0,0,x,y)=f(|y-x|\bmod 5)$, with $f(0)=\cP$,
$f(1)=\cN$, $f(2)=\cN$, $f(3)=\cP$, $f(4)=\cN$. Then (i) is immediate since
$|(y+t)-(x+t)|=|y-x|$.

For (ii), let $d=y-x$. Then $o(0,0,x+t,y)=f(|d-t|\bmod 5)$ and
$o(0,0,x,y)=f(|d|\bmod 5)$. We exhibit, for each residue of $t$ modulo $5$, a
value of $d$ with $f(|d-t|\bmod 5)\neq f(|d|\bmod 5)$:
\[
\begin{array}{c|ccccc}
t \bmod 5 & 0 & 1 & 2 & 3 & 4\\
\hline
d & t-3 & t-1 & t-2 & t+2 & t-4\\
|d-t| \bmod 5 & 3 & 1 & 2 & 2 & 4\\
f(|d-t|) & \cP & \cN & \cN & \cN & \cN\\
|d| \bmod 5 & 2 & 0 & 0 & 0 & 0\\
f(|d|) & \cN & \cP & \cP & \cP & \cP
\end{array}
\]
(Taking $d\ge 0$ here, which holds for $t\ge 4$; the finitely many smaller $t$
are checked directly.) In each case $f(|d-t|)\neq f(|d|)$, so taking $x$ large
with $y=x+d$ gives $o(0,0,x+t,y)\neq o(0,0,x,y)$, with $x,y$ arbitrarily large.
\end{proof}

\begin{corollary}\label{cor:03102infinite}
$\cQ(\texttt{0.3102})$ is infinite.
\end{corollary}

\begin{proof}
The slice $(0,0)$ is not bi--directionally eventually periodic
(Proposition~\ref{prop:03102slice}), so by the contrapositive of
Lemma~\ref{lem:bidirectional} the quotient is infinite.
\end{proof}

\begin{figure}[htbp]
\centering
\includegraphics[width=0.52\textwidth]{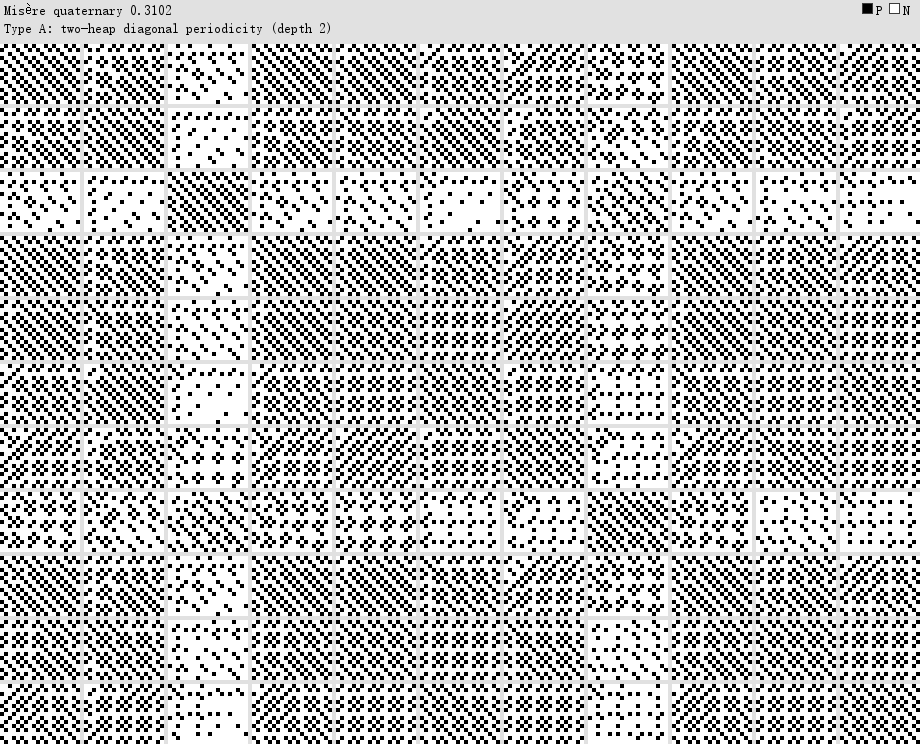}
\caption{The four--heap outcome of \texttt{0.3102} (range $20\times 15\times
10\times 10$, black $=\cP$, white $=\cN$).}
\end{figure}

\subsection{The infinitude of the mis\`ere quotients of \texttt{0.3122} and
\texttt{0.3123} (computer--assisted)}\label{sec:03122proven}

Unlike \texttt{0.3102}, the two--heap outcome of these two games is simply
periodic, and the non--periodicity appears only on \emph{four--heap slices}:
their infinite families are $(c,d,h)\not\equiv_{\cA}(c,d,h+1)$ ($c=d=2$ or $3$),
witnessed by the single heap $(h)$. The proof consists of two parts:
\begin{itemize}
\item \emph{finite verification}: the two reduction rules of slice--wise
algebraic periodicity are verified by exhaustive computation over a bounded
range (variable heaps $\le 150$, zero contradictions). This is the
\emph{computer--assisted} part---a consistency check of finite tables, too
large to spell out on paper, but fully reproducible;
\item \emph{hierarchical induction}: a purely mathematical argument---the
reduction decreases the total size, and the moves of each slice descend along a
partial order to lower slices, so strong induction on the total size extends the
finite verification to all positions.
\end{itemize}
Thus the mathematical skeleton of this proof agrees with the slice--aperiodicity
criterion for \texttt{0.3102}, except that the ``finite fact'' is obtained by
exhaustive verification (replacing the closed form $f(|y-x|\bmod 5)$ of
\texttt{0.3102})---which is exactly why this proof needs computer assistance and
is not fully written out on paper.

\begin{theorem}\label{thm:03122}
In \texttt{0.3122} (remove $\{1,2\}$, subtract $\{1,3,4\}$), for every $h\ge 0$,
\[
\om(2,2,h,h)=\cP,\qquad \om(2,2,h+1,h)=\cN .
\]
Consequently $\cQ(\texttt{0.3122})$ is infinite.
\end{theorem}

\begin{proof}
\emph{(1) Finite verification.} For a fixed prefix $(a,b)$ and variable heaps
$x,y$, write $x_1=\max(x,y)$, $x_2=\min(x,y)$. Exhaustive computation (variable
heaps $\le 150$, zero contradictions) confirms that the following slices satisfy
the two reduction rules of Definition~\ref{def:algperslice} (all with
$t_1=t_2=7$):
\[
\begin{array}{c|ccccccc}
(a,b) & (0,0) & (1,0) & (2,0) & (1,1) & (1,2) & (2,1) & (2,2)\\
\hline
(p_1,p_2) & (3,3) & (2,2) & (8,8) & (3,3) & (13,13) & (13,13) & (12,3)
\end{array}
\]

\emph{(2) Hierarchical induction.} The moves of these slices descend along the
partial order
\[
(0,0)\to(1,0)\to(2,0)\to(1,1)\to(1,2)\to(2,2),
\]
where $S\to S'$ means ``removing a fixed heap of $S$, or subtracting $1$ from a
fixed heap of $S$, lands in $S'$''. Since the reductions of (1) decrease the
total size, strong induction on the total size (running the slices in the above
order) establishes the reductions for all positions.

\emph{(3) The identities.} Apply the $(2,2)$--slice reduction to $(2,2,h,h)$:
for $h\ge 19$ the two--heap reduction $\om(2,2,h,h)=\om(2,2,h-7,h-7)$ applies
repeatedly, reducing to $(2,2,r,r)$ with $12\le r\le 18$ and $r\equiv h\pmod 7$.
The corner values (computed directly) satisfy $\om(2,2,r,r)=\cP$ and
$\om(2,2,r+1,r)=\cN$; the finitely many $h<19$ are checked directly.

\emph{(4) Infinite quotient.} From the identities,
$(2,2,h)\not\equiv_{\cA}(2,2,h+1)$ (witness $X=(h)$), and $\{(2,2,h)\}$ is
pairwise inequivalent.
\end{proof}

\begin{figure}[htbp]
\centering
\includegraphics[width=0.52\textwidth]{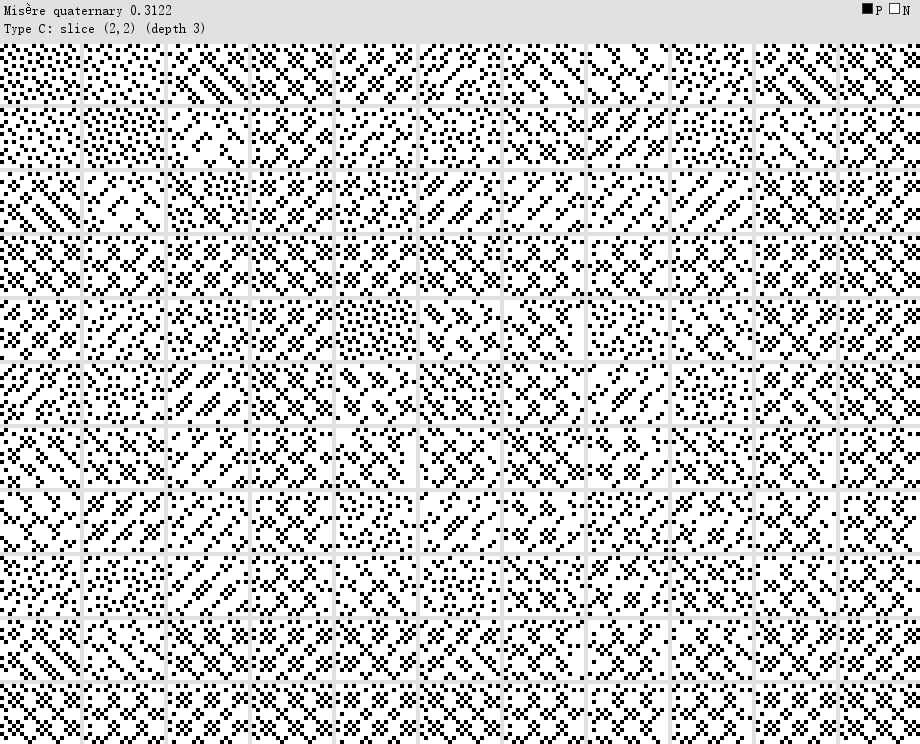}
\caption{The four--heap outcome of \texttt{0.3122} (black $=\cP$, white $=\cN$).}
\end{figure}

\begin{theorem}\label{thm:03123}
In \texttt{0.3123} (remove $\{1,2,4\}$, subtract $\{1,3,4\}$), for every
$h\ge 0$,
\[
\om(3,3,h,h)=\cP,\qquad \om(3,3,h+1,h)=\cN .
\]
Consequently $\cQ(\texttt{0.3123})$ is infinite.
\end{theorem}

\begin{proof}
The proof is completely parallel to Theorem~\ref{thm:03122}, except that the
hierarchical chain is longer. \emph{(1) Finite verification} confirms that the
following slices satisfy the two reduction rules of
Definition~\ref{def:algperslice} (all $t_1=t_2=7$, variable heaps $\le 150$, zero
contradictions):
\[
\begin{array}{c|cccccccccc}
(a,b) & (0,0) & (1,0) & (2,0) & (3,0) & (1,1) & (1,2) & (2,2) & (1,3) & (2,3) & (3,3)\\
\hline
(p_1,p_2) & (2,2) & (18,12) & (21,15) & (21,21) & (2,2) & (16,13) & (21,14) & (27,21) & (31,21) & (7,7)
\end{array}
\]
\emph{(2) Hierarchical induction} descends along the partial order
\[
(0,0)\to(1,0)\to(2,0)\to(3,0)\to(1,1)\to(1,2)\to(2,2)\to(1,3)\to(2,3)\to(3,3),
\]
and strong induction on the total size establishes the reductions for all
positions. \emph{(3) The identities}: applying the $(3,3)$--slice reduction to
$(3,3,h,h)$ reduces it to $(3,3,r,r)$ with $7\le r\le 13$ and
$r\equiv h\pmod 7$; the corner values satisfy $\om(3,3,r,r)=\cP$ and
$\om(3,3,r+1,r)=\cN$. \emph{(4) Infinite quotient}:
$(3,3,h)\not\equiv_{\cA}(3,3,h+1)$ (witness $X=(h)$).
\end{proof}

\begin{figure}[htbp]
\centering
\includegraphics[width=0.52\textwidth]{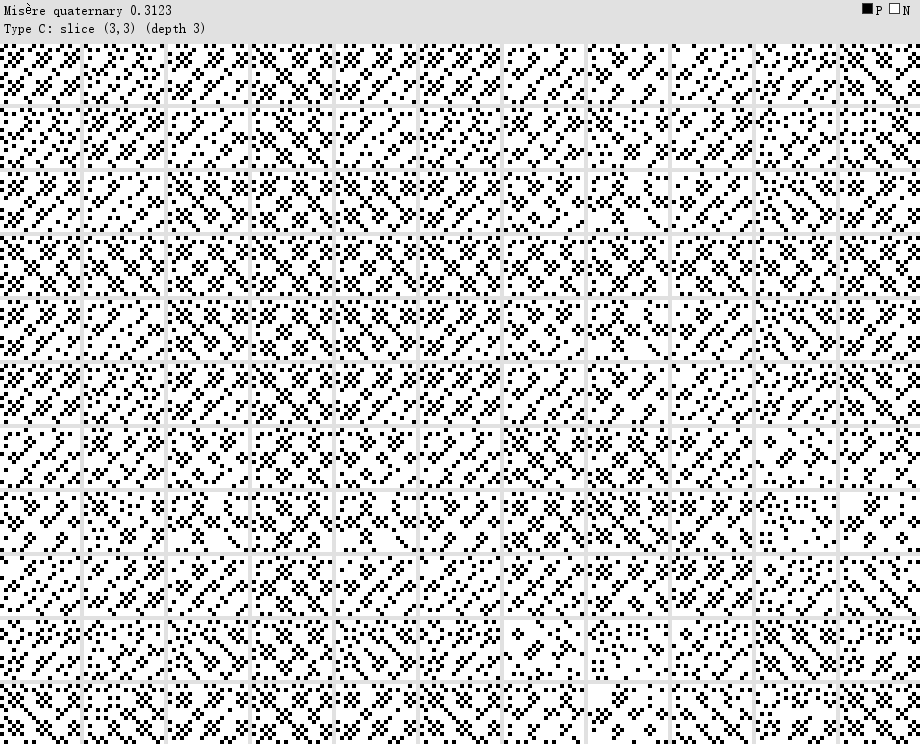}
\caption{The four--heap outcome of \texttt{0.3123} (black $=\cP$, white $=\cN$).}
\end{figure}

\section{Conjectures}\label{sec:conj}

In this section, the infinitude of the quotients of \texttt{0.3122} and
\texttt{0.3123} is already proved in Section~\ref{sec:proven} (as an application
of slice--wise algebraic periodicity and hierarchical induction); the remaining
results (including \texttt{0.3312} and the classification of the $8$--digit
games) are \emph{numerical only}, and are therefore stated as conjectures.

\subsection{Slice characterizations}\label{sec:slices}

A two--dimensional \emph{slice} $o(a,b,x,y)$ is the four--heap outcome with two
fixed heaps $a,b$ and two variable heaps $x,y$. The proofs of the infinitude of
the quotients of \texttt{0.3122}/\texttt{0.3123} (Section~\ref{sec:proven}) are
based on slices; in this subsection we give the regularities of the slices
themselves. The $98$--type characterization below is a numerically verified
conjecture (its diagonal and super--diagonal part is proved by
Theorem~\ref{thm:03122}); the rest are numerical observations.

\begin{conjecture}\label{conj:03122slice}
In \texttt{0.3122}, for all $x,y\ge 3$, the slice $o(2,2,x,y)$ depends only on
$\bigl(x\bmod 7,\ y\bmod 7,\ [\,|x-y|\le 4]\bigr)$, a table of
$7\times 7\times 2=98$ values.
\end{conjecture}

That is, the outcome is periodic in regions, with a different template near the
diagonal; this difference breaks the eventual periodicity, and the diagonal and
super--diagonal part is exactly the two identities of
Theorem~\ref{thm:03122}. The proof slice of \texttt{0.3123} is $(3,3)$
(Theorem~\ref{thm:03123}), whose slice--wise algebraic periodicity parameters
are $(p_1,t_1,p_2,t_2)=(7,7,7,7)$; no closed form similar to the $98$--type is
known.

For \texttt{0.3312}, no simple two--dimensional slice characterization was found:
its non--periodicity lies in the three--heap diagonal neighbourhood (all three
heaps variable) rather than in a two--dimensional slice, so the slice method of
Section~\ref{sec:proven} does not apply (see the open problems).

\subsection{Algebraic periodicity of the slices}\label{sec:algper}

We adapt Definition~\ref{def:algper} to a fixed slice, reducing only the two
\emph{variable} heaps.

\begin{definition}\label{def:algperslice}
For fixed $a,b$, the slice $o(a,b,x,y)$ is \emph{algebraically periodic} with
parameters $(p_1,t_1,p_2,t_2)$ if, writing $x_1=\max(x,y)$, $x_2=\min(x,y)$,
for all $x_1\ge p_1+t_1$:
\begin{enumerate}[label=(\roman*)]
\item if $x_2<p_2+t_2$, then $o(a,b,x_1,x_2)=o(a,b,x_1-t_1,x_2)$;
\item if $x_2\ge p_2+t_2$, then $o(a,b,x_1,x_2)=o(a,b,x_1-t_2,x_2-t_2)$.
\end{enumerate}
\end{definition}

Thus the \emph{bulk} (both variable heaps large) is diagonal--periodic with
period $t_2$, and the \emph{edge} (one variable heap small) is axis--periodic
with period $t_1$.

\begin{conjecture}\label{conj:algper}
For each of the four games \texttt{0.3102}, \texttt{0.3122}, \texttt{0.3123},
\texttt{0.3312}, taking $t_1=t_2=$ its single--heap period (respectively
$5,7,7,3$), the following slices are \emph{exactly} algebraically periodic (zero
contradictions over the tested range):
\begin{itemize}
\item the $2$--heap slice $(0,0)$;
\item every $3$--heap slice $(a,0)$ and $(0,b)$;
\item the diagonal $4$--heap slices $(a,a)$.
\end{itemize}
\end{conjecture}

\begin{conjecture}\label{conj:foldfail}
Folding the fixed heaps by the period fails: in general
$o(a+p,b,x,y)\neq o(a,b,x,y)$ (contradictions $561$--$1935$ over the tested
range). Hence the fixed heaps cannot be reduced modulo the period, consistently
with the infinite quotient.
\end{conjecture}

The clean stratification is:

\begin{center}
\begin{tabular}{lc}
\toprule
slice & algebraic periodicity\\
\midrule
$2$--heap $(0,0)$ & exact\\
$3$--heap $(a,0),(0,b)$ & exact\\
$4$--heap diagonal $(a,a)$ & exact\\
$4$--heap off--diagonal $(a\neq b)$ & exact (large pre--period)\\
\bottomrule
\end{tabular}
\end{center}

\begin{remark}
In fact, taking the pre--period $p_1$ and the boundary $p_2$ sufficiently large,
\emph{all} slices (including the off--diagonal ones) of the four games are
exactly algebraically periodic (zero contradictions over the tested range): for
example \texttt{0.3123} admits $(2,3)\mapsto(31,21)$ and
$(3,4)\mapsto(31,21)$, and \texttt{0.3102} admits $(1,2)\mapsto(46,12)$ and
$(2,3)\mapsto(46,12)$. Hence what distinguishes the games is \emph{not} whether
some slice is algebraically periodic (they all are), but \emph{which} slice's
diagonal fails to be eventually periodic: $(0,0)$ for \texttt{0.3102},
$(2,2)/(3,3)$ for \texttt{0.3122}/\texttt{0.3123}; for \texttt{0.3312} the
three--heap diagonal neighbourhood is aperiodic, so the infinitude of its
quotient still requires a new idea.
\end{remark}

\subsection{Wildness depth}\label{sec:depth}

The \emph{wildness depth} is the minimum number of heaps needed to make an
explicit infinite family $G_h\not\equiv_{\cA}G_{h+p}$ pairwise inequivalent,
namely the number of heaps of the witness $X_h$ together with the family member
$G_h$---when $G_h$ is a single heap, this counts the $1$ heap of $G_h$ itself
plus the heaps of the witness; when $G_h$ is a multi--heap, the witness is
already contained in $G_h$, and the depth is the number of heaps of $G_h$. For
example, the family of \texttt{0.3102} is the single heap $(h)$ with the single
--heap witness $(h+2)$, depth $1+1=2$; the family of \texttt{0.3122} is
$(2,2,h)$ ($3$ heaps) with witness $(h)$ contained in it, depth $3$; the family
of \texttt{0.3312} is the single heap $(h)$ with the two--heap witness $(h,h)$,
depth $1+2=3$.

\begin{center}
\begin{tabular}{ccc}
\toprule
Game & wildness depth & witness type\\
\midrule
\texttt{0.3102} & $2$ & single heap $(h+2)$\\
\texttt{0.3122} & $3$ & single heap $(h)$\\
\texttt{0.3123} & $3$ & single heap $(h)$\\
\texttt{0.3312} & $3$ & two heaps $(h,h)$\\
\bottomrule
\end{tabular}
\end{center}

\begin{figure}[htbp]
\centering
\includegraphics[width=0.52\textwidth]{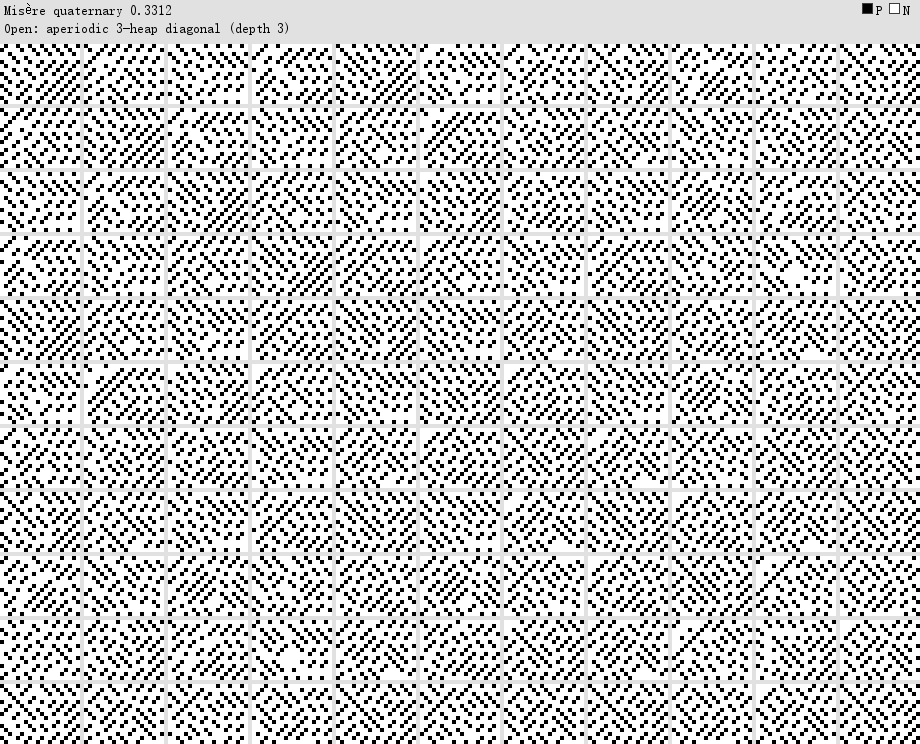}
\caption{The four--heap outcome of \texttt{0.3312} (black $=\cP$, white $=\cN$).}
\end{figure}

\begin{conjecture}\label{conj:depth}
The wildness depth of each game is as above. Moreover (all verified
numerically) no simpler finite description exists: (i) the three--heap outcome
is not eventually periodic in any coordinate (periods $3,4,5,7$ all fail);
(ii) the mis\`ere outcome is neither a function of the normal--play Grundy
profile nor of the normal nim--sum (so the normal kernel is infinite);
(iii) there is no translation structure $(h+p)\equiv_{\cA}(h)+(p)$ for any
$p\le 8$.
\end{conjecture}

\section{Wild quaternary games with longer codes}\label{sec:longer}

We extend the search to \emph{all} quaternary codes of length at most $8$,
$\Gamma=0.d_1d_2\cdots d_8$ with $d_1,d_2,\dots,d_8\in\{0,1,2,3\}$ (no restriction
on the first digit), with trailing zeros stripped and the trivial all--zero code
excluded ($65535$ distinct normalized codes, including the $1$-- through
$8$--digit games). For each game we detect whether the single--heap images
$\Phi((n))$ are eventually periodic: we take a set $W$ of witness positions, call
the tuple $(\om(G+X))_{X\in W}$ the \emph{signature} of a position $G$, and check
whether the signatures for $n=1,\dots,50$ are eventually periodic. For the
aperiodic ones we then search for an explicit infinite family
$G_h\not\equiv_{\cA}G_{h+p}$.

\emph{Caveat.} The content of this section is \emph{numerical only}: all the
``wild''/``tame'' classifications, the explicit families, and the counts below
are obtained by finite computation over the stated ranges; no rigorous proof is
currently known, so absolute accuracy is \emph{not} guaranteed. In particular,
for a game labelled ``confirmed wild'', the pairwise inequivalence of its
explicit family is verified only for $h\le 40$; the \emph{reduction} (an
infinite family $\Rightarrow$ infinite quotient) is a rigorous Myhill--Nerode
argument, but the family itself is only numerically verified, and its validity
for all $h$ is not proved.

\begin{remark}[A false--periodicity pitfall]\label{rem:falseper}
The ``first repeat'' periodicity test must be applied with care: if it accepts a
period based on only a short tail of the signature sequence, it can misclassify
a wild game (whose signature keeps producing new values as $n$ grows) as tame.
We therefore require the period to hold over a tail of at least $15$ values.
Without this, e.g.\ \texttt{0.3312} is spuriously declared to have period $3$,
and $115$ wild games are missed.
\end{remark}

\begin{conjecture}\label{conj:eightdigit}
Among the $65535$ quaternary games of length at most $8$ (with no restriction on
the first digit), $15279$ are \emph{numerically confirmed wild}: each admits an
explicit infinite family $(h)\not\equiv_{\cA}(h+p)$ (verified for $h\le 40$),
whose witness varies with $h$ (of the form $(h+a)$ or $(h,h)$). A further $735$
games are \emph{candidates} (aperiodic single--heap images, but no explicit
family was found). The $7$-- and $8$--digit games alone contribute $14805$
confirmed wild and $723$ candidates, i.e.\ roughly $24\%$ of all $7$--/
$8$--digit codes are numerically classified as wild. These numbers and
classifications are subject to the caveat above.
\end{conjecture}

\begin{remark}[Provable wildness classification]\label{rem:provable}
Reclassifying the $14805$ numerically confirmed wild games of length $7$--$8$
by the rigorous methods of this paper, we obtain:
\begin{itemize}
\item \textbf{Type A} ($0.3102$--type; wildness depth $2$): the two--heap
outcome is diagonally periodic, $o(x,y)=f(|y-x|\bmod p)$. Proved by the
slice--aperiodicity criterion. \textbf{$2$ games.}
\item \textbf{Type B} ($0.331012$--type, ``maximally wild''; depth $2$):
$o(h,h)=\cP$ and $o(h+p,h)=\cN$. Proved by the clean induction of
Theorem~\ref{thm:0331012}. \textbf{$409$ games.}
\item \textbf{Type C} ($0.3122$--type; depth $3$): some slice $(c,d)$
($c,d\ge 1$) satisfies $o(c,d,h,h)\neq o(c,d,h+1,h)$, giving the family
$(c,d,h)\not\equiv_{\cA}(c,d,h+1)$ (witness $(h)$). Proved by slice--wise
algebraic periodicity and hierarchical induction. \textbf{$5582$ games.}
\item \textbf{Numerical only}: the remaining \textbf{$8812$ games}, whose
wildness is currently established only by an explicit family verified for
$h\le 40$.
\end{itemize}
Thus $5993$ of the $14805$ wild games are \emph{provably} wild.
\end{remark}

The wild games are distributed across all four possible first digits
$d_1\in\{0,1,2,3\}$ (counts for the $7$-- and $8$--digit games, numerically
confirmed):
\begin{itemize}
\item $d_1=0$ (prefix $\texttt{0.0}$, no move on a single heap): $4027$ wild
games, e.g.\ $\texttt{0.0122}$;
\item $d_1=1$ (prefix $\texttt{0.1}$, ``remove $1$'' only): $5184$ wild games;
\item $d_1=2$ (prefix $\texttt{0.2}$, ``subtract $1$'' only): $2890$ wild games;
\item $d_1=3$ (prefix $\texttt{0.3}$, ``remove or subtract $1$''): $2704$ wild
games, containing the four known games $\texttt{0.3102},\texttt{0.3122},
\texttt{0.3123},\texttt{0.3312}$ and all their longer extensions.
\end{itemize}
Thus the four known games are only the ``tip'' of a much larger wild family: the
``$d_1=1$'' family and, more surprisingly, the ``$d_1=0$'' family (with no move
on a single heap) are even more numerous than the $d_1=3$ extensions.

Many of these are ``maximally wild'': the single heaps are pairwise
inequivalent ($(h)\not\equiv_{\cA}(h+1)$), witnessed by $(h)$. For such games the
following clean induction gives a rigorous proof.

\begin{theorem}\label{thm:0331012}
In the game $\texttt{0.331012}$ (remove $\{1,2,3,5\}$, subtract $\{1,2,6\}$),
for every $h\ge 3$,
\[
\om(h,h)=\cP,\qquad \om(h+1,h)=\cN .
\]
\end{theorem}

\begin{proof}
By strong induction on $h$. The base cases $h=3,\dots,8$ are checked directly
(twelve values, all satisfying $\om(h,h)=\cP$ and $\om(h+1,h)=\cN$). For
$h\ge 9$ we have $h\notin\{1,2,3,5\}$ and $h>6$, so the options of $(h,h)$ are
exactly $(h-1,h)$, $(h-2,h)$, $(h-6,h)$ (subtracting $1$, $2$, $6$ from one
heap; no removal). Each is an $\cN$--position:
\begin{itemize}
\item $(h-1,h)$ has the option $(h-1,h-1)$ (subtract $1$ from the larger heap),
which is $\cP$ by induction at $h-1$;
\item $(h-2,h)$ has the option $(h-2,h-2)$ (subtract $2$), $\cP$ by induction
at $h-2$;
\item $(h-6,h)$ has the option $(h-6,h-6)$ (subtract $6$), $\cP$ by induction
at $h-6$.
\end{itemize}
Hence all options of $(h,h)$ are $\cN$, so $\om(h,h)=\cP$. Then $(h+1,h)$ has
the option $(h,h)$ (subtract $1$ from the larger heap), which is $\cP$, so
$\om(h+1,h)=\cN$. This closes the induction.
\end{proof}

\begin{corollary}
$(h)\not\equiv_{\cA}(h+1)$ (witness $X=(h)$, as $\om(h,h)=\cP\neq\cN=
\om(h+1,h)$), hence the mis\`ere quotient of $\texttt{0.331012}$ is infinite.
The same proof applies verbatim to $\texttt{0.331002}$ and $\texttt{0.331211}$
(whose subtract sets are $\{1,2,6\}$ and $\{1,2,4\}$).
\end{corollary}

\begin{figure}[htbp]
\centering
\includegraphics[width=0.52\textwidth]{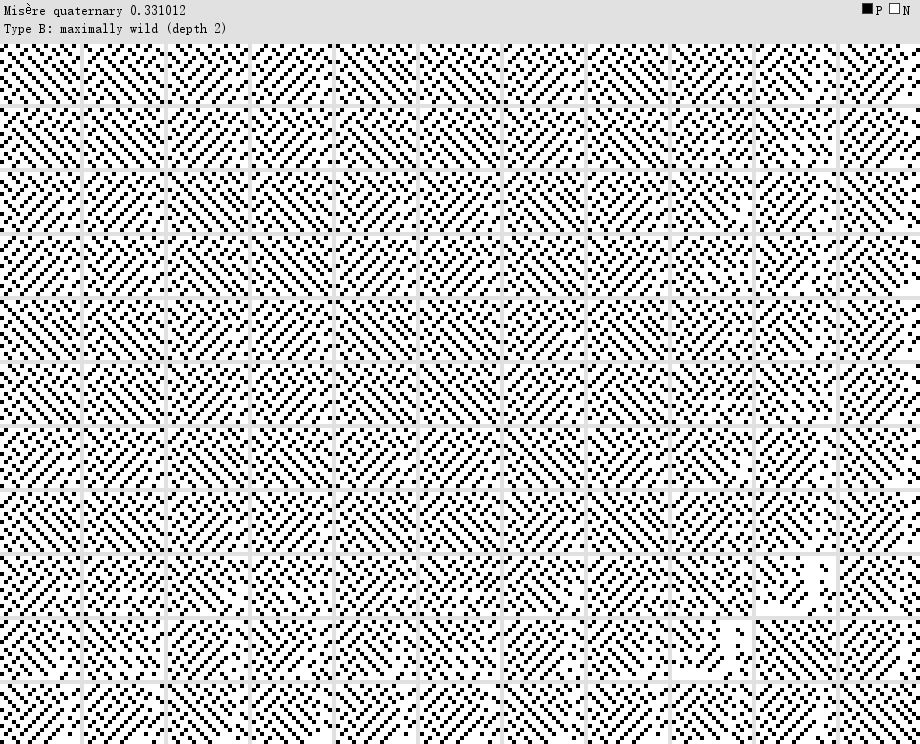}
\caption{The four--heap outcome of \texttt{0.331012} (black $=\cP$, white
$=\cN$).}
\end{figure}

\section{Open problems}\label{sec:open}

\begin{problem}
Prove that the quotient of \texttt{0.3312} is infinite: its identities are
$\om(h,h,h)=\cN\iff h\equiv0\pmod 3$ and $\om(h+4,h,h)=1-\om(h,h,h)$ (verified
for $h\ge 19$). But the full three--heap algebraic periodicity fails: even with
the optimal parameters $(p_1,t_1,p_2,t_2)=(29,3,29,3)$ there are $903$
contradictions, e.g.\ $\om(33,33,29)=\cP\neq\cN=\om(30,30,29)$, while the
diagonal $\om(h,h,h)=[h\equiv0\bmod 3]$ has period $3$. Thus the diagonal
neighbourhood is aperiodic, the slice--wise hierarchical induction does not
close, and a new idea is needed.
\end{problem}

\begin{problem}
Characterize the bulk of the slice $o(1,2,x,y)$ of \texttt{0.3123}: it admits no
description by $\bigl(x\bmod 7,y\bmod 7,[\,|x-y|\le t]\bigr)$. What is the
correct invariant?
\end{problem}

\begin{problem}
Determine the precise automaton class that recognizes the outcomes of these
games: finite automata do not suffice (\cite{gamedfa}), since the outcome is not
even eventually periodic in any coordinate; does a \emph{pushdown automaton}
(an automaton with an unbounded stack) suffice? Does a ``difference + residue''
representation (taking the differences of the heaps, rather than their absolute
values, as coordinates) make the outcome recognizable by a finite automaton?
\end{problem}

\begin{problem}
Determine whether the ``weak algebraic periodicity'' of \cite{gamedfa} for
$\mathbf{0.145}$ and $\mathbf{0.54}$ is genuine, or likewise an artefact of a
too small pre--period---given that, for the four quaternary games, \emph{all}
slices are exactly algebraically periodic once the pre--period is large.
\end{problem}

\begin{problem}
Explain why \emph{all} slices (including the off--diagonal ones) of the four
games are exactly algebraically periodic once the pre--period is large, and why
what really distinguishes the games is \emph{which} slice's diagonal fails to be
eventually periodic ($(0,0)$ for \texttt{0.3102}, $(2,2)/(3,3)$ for
\texttt{0.3122}/\texttt{0.3123}, and the three--heap diagonal neighbourhood for
\texttt{0.3312}).
\end{problem}

\vspace{1em}
\noindent\textbf{Acknowledgments.} I thank Martin Hans Keskküla (University of
Tartu) for pointing out the counterexample $(2,2,9,4)$ and $(2,2,2,11)$ to the
$98$--type characterization, which led to the lower bound $x,y\ge 3$ in
Conjecture~\ref{conj:03122slice}.

\vspace{1em}
\noindent\textbf{AI usage statement.} An AI tool (DeepSeek) was used during the
preparation of this paper to assist with drafting and proofreading parts of the
text, with writing the numerical--computation and verification programs, and
with literature search. The correctness of all mathematical proofs and numerical
results has been independently checked by the author, who is responsible for the
entire content.

\bibliographystyle{plain}
\bibliography{refs}

@article{PS,
  author        = {Plambeck, T. E. and Siegel, A. N.},
  title         = {Mis\`ere quotients for impartial games},
  journal       = {Journal of Combinatorial Theory, Series A},
  volume        = {115},
  number        = {4},
  pages         = {593--622},
  year          = {2008},
  eprint        = {math/0609825},
  archivePrefix = {arXiv},
  url           = {https://arxiv.org/abs/math/0609825}
}

@misc{miseregames,
  author       = {Siegel, A. N.},
  title        = {Octal games of the form 0.3xxx},
  howpublished = {miseregames.org},
  url          = {http://miseregames.org/solutions/summary.php?name=0.3xxx}
}

@misc{gamedfa,
  author        = {Liang, Kai},
  title         = {A Finite Automaton Approach to Combinatorial Games},
  year          = {2026},
  eprint        = {2608.13273},
  archivePrefix = {arXiv},
  primaryClass  = {math.CO},
  url           = {https://arxiv.org/abs/2608.13273}
}

@article{MisQuot2006,
  author  = {Plambeck, T. E.},
  title   = {Taming the wild in impartial combinatorial games},
  journal = {Integers},
  volume  = {5},
  number  = {1},
  pages   = {G5},
  year    = {2005}
}

\end{document}